\documentclass[a4paper,12pt]{amsart}
\usepackage{amssymb, amsfonts, amsthm}
\usepackage{ifthen}
\usepackage{graphicx}
\nonstopmode \numberwithin{equation}{section}
\usepackage[usenames]{color}
\newtheorem{thm}{Theorem}
\newtheorem{cor}{Corollary}
\newtheorem{lem}{Lemma}

\newtheorem{rem}{Remark}
\newtheorem{rems}[equation]{Remarks}
\theoremstyle{definition}
\newtheorem{defin}{Definition}
\newtheorem{examp}[equation]{Example}
\newtheorem{prob}[equation]{Problem}
\newtheorem{ques}[equation]{Question}
\newtheorem{op}{Open Problem}

\newtheorem{conj}[equation]{Conjecture}
\newtheorem{deter}[equation]{Determination}
\newtheorem{case}{Case}[section]
\newtheorem{subcase}[equation]{Subcase}
\newtheorem{claim}{Claim}[section]
\newtheorem{subclaim}{Subclaim}

\newcounter {own}
\def\theown {\thesection       .\arabic{own}}

\newenvironment{pf}[1][]{%
 \vskip 3mm
 \noindent
 \ifthenelse{\equal{#1}{}}%
  {{\bf Proof. }}%
  {{\bf #1.} }%
 }%
{\qed\bigskip}

\newcounter{alphabet}
\newcounter{tmp}

\makeatletter
\newcommand{\Ref}[1]{\@ifundefined{r@#1}{}{\setcounter{tmp}{\ref{#1}}\Alph{tmp}}}
\makeatother

\newenvironment{Lem}[1][]{\refstepcounter{alphabet}%
\bigskip%
\noindent%
{\bf Lemma \Alph{alphabet}}%
{\bf .} \itshape}{\vskip 8pt}

\newcommand{\IR}{{\mathbb R}}

\newcommand{\diam}{{\operatorname{diam}}}

\def\be{\begin{equation}}
\def\ee{\end{equation}}

\newcommand{\bee}{\begin{enumerate}}
\newcommand{\eee}{\end{enumerate}}

\newcommand{\blem}{\begin{lem}}
\newcommand{\elem}{\end{lem}}
\newcommand{\bthm}{\begin{thm}}
\newcommand{\ethm}{\end{thm}}
\newcommand{\bcor}{\begin{cor}}
\newcommand{\ecor}{\end{cor}}
\newcommand{\beg}{\begin{examp}}
\newcommand{\eeg}{\end{examp}}
\newcommand{\begs}{\begin{examples}}
\newcommand{\eegs}{\end{examples}}
\newcommand{\bdefe}{\begin{defin}}
\newcommand{\edefe}{\end{defin}}
\newcommand{\bprob}{\begin{prob}}
\newcommand{\eprob}{\end{prob}}
\newcommand{\bques}{\begin{ques}}
\newcommand{\eques}{\end{ques}}
\newcommand{\bei}{\begin{itemize}}
\newcommand{\eei}{\end{itemize}}

\newcommand{\bde}{\begin{deter}}
\newcommand{\ede}{\end{deter}}
\newcommand{\bca}{\begin{case}}
\newcommand{\eca}{\end{case}}
\newcommand{\bsca}{\begin{subcase}}
\newcommand{\esca}{\end{subcase}}
\newcommand{\bcl}{\begin{claim}}
\newcommand{\ecl}{\end{claim}}
\newcommand{\bscl}{\begin{subclaim}}
\newcommand{\escl}{\end{subclaim}}
\newcommand{\bcon}{\begin{conj}}
\newcommand{\econ}{\end{conj}}
\newcommand{\bcons}{\begin{conjs}}
\newcommand{\econs}{\end{conjs}}
\newcommand{\bprop}{\begin{propo}}
\newcommand{\eprop}{\end{propo}}
\newcommand{\br}{\begin{rem}}
\newcommand{\er}{\end{rem}}
\newcommand{\brs}{\begin{rems}}
\newcommand{\ers}{\end{rems}}
\newcommand{\bo}{\begin{obser}}
\newcommand{\eo}{\end{obser}}
\newcommand{\bos}{\begin{obsers}}
\newcommand{\eos}{\end{obsers}}
\newcommand{\bpf}{\begin{pf}}
\newcommand{\epf}{\end{pf}}
\newcommand{\ba}{\begin{array}}
\newcommand{\ea}{\end{array}}
\newcommand{\beq}{\begin{eqnarray}}
\newcommand{\beqq}{\begin{eqnarray*}}
\newcommand{\eeq}{\end{eqnarray}}
\newcommand{\eeqq}{\end{eqnarray*}}

\newcommand{\bop}{\begin{op}}
\newcommand{\eop}{\end{op}}

\newtheorem{pfofThm1.5}[equation]{}

\newcounter{minutes}
\divide\time by 60
\newcounter{hours}
\multiply\time by 60 \addtocounter{minutes}{-\time}

\begin{document}

\bibliographystyle{amsplain}

\title{Weakly quasisymmetric maps and uniform spaces}

\def\thefootnote{}
\footnotetext{ \texttt{\tiny File:~\jobname .tex,
          printed: \number\year-\number\month-\number\day,
          \thehours.\ifnum\theminutes<10{0}\fi\theminutes}
} \makeatletter\def\thefootnote{\@arabic\c@footnote}\makeatother

\author{Yaxiang  Li}
\address{Yaxiang Li,  College of Science,
Central South University of
Forestry and Technology, Changsha,  Hunan 410004, People's Republic
of China} \email{yaxiangli@163.com}

\author[]{Matti Vuorinen}
\address{Matti Vuorinen, Department of Mathematics and Statistics, University of Turku,
FIN-20014 Turku, Finland}
\email{vuorinen@utu.fi}


\author{Qingshan Zhou${}^{\mathbf{*}}$}
\address{Qingshan Zhou, school of mathematics and big data, foshan university,  Foshan, Guangdong 528000, People's Republic
of China} \email{q476308142@qq.com}

\date{}
\keywords{Subinvariance, uniform domain, weak quasisymmetry, short arc,quasihyperbolic geodesic.
\\
${}^{\mathbf{*}}$ Corresponding author}

\begin{abstract}
Suppose that $X$ and $Y$ are quasiconvex and complete metric spaces, that $G\subset X$ and $G'\subset Y$ are
domains, and that $f: G\to G'$ is a homeomorphism. In this paper, we first give some basic properties of short arcs, and then we show that: if $f$ is a weakly quasisymmetric mapping and $G'$ is
a quasiconvex domain, then the image $f(D)$ of every uniform subdomain $D$ in $G$  is uniform. As an application, we get that if $f$ is a weakly quasisymmetric mapping and $G'$ is
an uniform domain, then the images of the short arcs in $G$ under $f$ are uniform arcs in the sense of diameter.
\end{abstract}

\thanks{The research was partly supported by NNSF of
China (No. 11601529, No. 11671127) and NSF of Hunan Province (No. 2015JJ3171).}

\maketitle\pagestyle{myheadings} \markboth{Yaxiang Li, Matti Vuorinen and Qingshan Zhou}{Weakly quasisymmetric maps preserves uniform domains in metric spaces}

\section{Introduction and main results}\label{sec-1}

The quasihyperbolic metric (briefly, QH metric) was introduced by Gehring and
his students Palka and Osgood in the 1970's \cite{Geo, GP} in the setting of Euclidean spaces
${\mathbb R}^n$ $(n\ge 2).$ Since its first appearance, the quasihyperbolic metric has become an important tool in the geometric function theory of Euclidean spaces, especially, in the study of quasiconformal and quasisymmetric mappings. Uniform domains in Euclidean spaces were introduced independently by Jones \cite{Jo80} and Martio and Sarvas \cite{MS}. Recently, Bonk, Heinonen and Koskela introduced uniform metric spaces in \cite{BHK} and demonstrated a one-to-one (conformal) correspondence between this class of spaces and geodesic hyperbolic spaces in the sense of Gromov. After its appearance, uniformity has played a significant role in related studies; see \cite{BH}, \cite{Her04}, \cite{Her06}, \cite{KL}, \cite{KLM14}, \cite{La11} and references therein.

The class of quasisymmetric mappings on the real axis was first introduced by Beurling and
Ahlfors \cite{BA}, who found a way to obtain a quasiconformal extension of a quasisymmetric self-mapping of the
real axis to a self-mapping of the upper half-plane. This idea
was later generalized by Tukia and V\"ais\"al\"a, who studied quasisymmetric
mappings between metric spaces \cite{TV}. In
1998, Heinonen and Koskela \cite{HK} proved a remarkable result, showing that the concepts of quasiconformality and
quasisymmetry are quantitatively equivalent in a large class of metric spaces, which
includes Euclidean space.  Also, V\"ais\"al\"a proved the quantitative equivalence between free quasiconformality and
quasisymmetry of homeomorphisms between two Banach spaces, see \cite[Theorem 7.15]{Vai8}.  Against this background, it is not surprising that the study of quasisymmetry in metric spaces has recently attracted significant attention \cite{BM,hprw,HL,Tys, WZ}.

The main tools in the study of quasiconformal mappings and uniform spaces are volume integrals (associated doubling or Ahlfors regular measure), conformal modulus, Whitney decomposition and quasihyperbolic metric. The main goal of this paper is to study the subinvariance of uniform domains in general metric spaces under weakly quasisymmetric mappings by means of the quasihyperbolic metric and metric geometry. We start by recalling some basic definitions.
Throughout this paper, we always assume that $X$ and $Y$ are metric spaces and we do not assume local compactness. We follow the notation and terminology of  \cite{HK-1, HK, HL, Tys, Vai8}.
%
\noindent Here and in what follows, we always use $|x-y|$ to denote the distance between $x$ and $y$.

\bdefe\label{japan-31} A homeomorphism $f$ from $X$ to $Y$ is said to be
\begin{enumerate}
\item $\eta$-{\it quasisymmetric} if there is a homeomorphism $\eta : [0,\infty) \to [0,\infty)$ such that
$$ |x-a|\leq t|x-b|\;\; \mbox{implies}\;\;   |f(x)-f(a)| \leq \eta(t)|f(x)-f(b)|$$
for each $t\geq 0$ and for each triple $x,$ $a$, $b$ of points in $X$;

\item {\it weakly $H$-quasisymmetric} if
$$ |x-a|\leq |x-b|\;\;  \mbox{ implies}\;\;   |f(x)-f(a)| \leq H|f(x)-f(b)|$$
for each triple $x$, $a$, $b$ of points in $X$.
\end{enumerate}\edefe

\br\label{japan-32}
The $\eta$-quasisymmetry implies the weak $H$-quasisymmetry with $H=\eta(1)$. Obviously, $\eta(1)\geq 1$. In general, the converse is not true (cf. \cite[Thm. $8.5$]{Vai8}). See also \cite{hprw} for some related results.
\er

In his 1961 work on
elasticity theory, John  \cite{John} introduced the class of domains satisfying the twisted interior cone condition . These domains were first called John domains by Martio and Sarvas
in \cite{MS}. In the same paper, Martio and Sarvas also discussed another class of domains which are the uniform domains. Their main motivation
 for studying these domains was in showing global injectivity properties for locally injective mappings. Since then, many other characterizations of uniform and John domains have been established, see
\cite{FW,  Geo, Martio-80, Vai6, Vai4, Vai8}, and the importance of these classes of domains in function theory is well
documented (see e.g. \cite{FW, GH, Vai2}). Moreover, John and uniform domains in
$\mathbb{R}^n$ enjoy numerous geometric and function theoretic
properties that are useful in many other fields of modern mathematical analysis as well (see e.g.
\cite{Jo80, Yli, Vai2}, and references therein).

We recall the definition of uniform domains  following closely the notation
and terminology of \cite{TV, Vai2, Vai, Vai6-0, Vai6}  and \cite{Martio-80}.

\bdefe \label{def1.3} A domain $G$ in $X$ is called $b$-{\it
uniform} provided there exists a constant $b$
with the property that each pair of points $x$, $y$ in $G$ can
be joined by a rectifiable arc $\gamma$ in $G$ satisfying
 \bee
\item $\min\{\ell(\gamma[x,z]),\ell(\gamma[z,y])\}\leq b\,\delta_{G}(z)$ for all $z\in \gamma$, and
\item $\ell(\gamma)\leq b\,|x-y|$,
\eee
\noindent where $\ell(\gamma)$ denotes the length of $\gamma$,
$\gamma[x,z]$ the part of $\gamma$ between $x$ and $z$.
In  this case, $\gamma$ is said to be a {\it double $b$-cone arc}.
If the condition $(1)$ is satisfied, not necessarily $(2)$, then $G$ is said to be a {\it $b$-John domain}. In this case, the arc $\gamma$ is called a {\it $b$-cone arc}.
\edefe

It follows from \cite[Rem. on p. 121]{FHM} and \cite[Thm. 5.6]{Vai2} that uniform domains are subinvariant with respect to quasiconformal mappings in $\IR^n$ ($n\geq 2$). By this, we mean that if $f:\; G\to G'$ is a $K$-quasiconformal mapping, where $G$ and $G'$ are domains in $\IR^n$, and if $G'$ is $c$-uniform,
then $D'=f(D)$ is $c'$-uniform
 for every $c$-uniform subdomain
 $D \subset G$,
 where $c'=c'(c, K, n)$ which means that the constant $c'$ depends only on the coefficient  $c$ of the uniformity of $D$, the coefficient $K$ of quasiconformality of $f$ and the dimension $n$ of the Euclidean space $\IR^n$. See \cite{BHX, GM, hlpw, HVW, Vai2, Vai8, Xie} for other relevant discussions. We note that a domain $G$ is uniform implies that $G$ is a John domain and is quasiconvex. So it is natural to ask whether it is possible to weaken the assumption ``$G'$ is uniform" to ``$G'$ is a John domain" or ``$G'$ is quasiconvex".

In fact, we observe from \cite[Thm. 1]{hlpw} and \cite[Prop. 7.12]{BHK} that the following   holds:
{\it Suppose that $G$ and $G'$ are bounded subdomains in $\mathbb{R}^n$ and that $f:G\to G'$ is a $K$-quasiconformal mapping. If $G'$ is a John domain and  $G_1$ is a subdomain of $G$ which is inner  uniform, then its image $G_1'=f(G_1)$ is inner uniform also.}
We remark that this result is not valid for uniform subdomain $G_1$ of $G$, that is, $G_1'$ maybe not uniform. For example, $G'=\mathbb{B}^2\setminus [0,1]$ is a conformal image of $G=\mathbb{B}^2$, and we observe that $G'$ is a John domain and $G_1=G\setminus\{0\}$ is uniform, but $G_1'$ is obviously not an uniform domain, because it is of the form $G'$ minus one point.  However, if we replace the assumption ``$G'$ is an $a$-John domain" to ``$G'$ is quasiconvex", then we get the following result.

\begin{thm}\label{thm1.1}
Suppose that $X$ and $Y$ are quasiconvex and complete metric spaces, that $G\varsubsetneq X$ is a domain,
$G'\varsubsetneq Y$ is a quasiconvex domain, and that  $f: G\to G'$ is weakly quasisymmetric mapping.
For each subdomain $D$ of $G$, if $D$ is uniform, then $D'=f(D)$ is uniform, where the coefficient of uniformity of $D'$ depends
only on the given data of $X$, $Y$, $G$, $G'$, $D$, and $f$.
\end{thm}

 Here and in what follows,
the phrase ``the given data  of $X$, $Y$, $G$, $G'$, $D$, and $f$" means the data which depends on the given constants which are the coefficients of quasiconvexity of $X$, $Y$ and $G'$, the coefficient of uniformity of $G$
and  the coefficient of weak quasisymmetry of $f$.

\br It is worth mentioning that in Theorem \ref{thm1.1}, the domain $G'$ is not required to be ``uniform", but only to be ``quasiconvex" $($From the definitions
in Section \ref{sec-2}, we easily see that uniformity implies quasiconvexity$)$.  If $X=Y=\mathbb{R}^n$, then $f$ is $\eta$-quasisymmetric with $\eta=\eta(n,H)$, see \cite[Thm. $2.9$]{Vai0}. Since quasisymmetric maps preserve uniform domains, the assertion follows. But we remind the reader that our result is independent of the dimension in this case.

\er

As an application of our method, we discuss the distortion property for quasihyperbolic short arcs because in general the quasihyperbolic geodesics may not exist. Actually we establish an analog of Pommerenke's theorem for length and diameter distortion of hyperbolic geodesics under conformal mappings from the unit disk $D$ onto a plane domain, see \cite[Cor. $4.18$, Thm. $4.20$]{Po}. In higher dimension (that is, $\mathbb{R}^n$), Heinonen and N\"{a}kki found that the quasiconformal image of quasihyperbolic geodesics minimize Euclidean curve-diameter, see \cite[Thm. $6.1$]{HN}. We also get the diameter uniformity for the image of quasihyperbolic short arcs under weakly quasisymmetric mappings which embed into a uniform space, which we state as follows.

\begin{thm}\label{thm-3}
Suppose that $X$ and $Y$ are $c$-quasiconvex and complete metric spaces, and that $f:G\to G'$ is weakly $H$-quasisymmetric between two domains $G\subsetneq X$ and $G'\subsetneq Y$. If $G'$ is $b$-uniform, then for any $\varepsilon$-short arc $\gamma$ in $G$ with endpoints $x$ and $y$, $0<\varepsilon<\min\{1,\frac{k_G(x,y)}{6}\}$, we have
\begin{enumerate}
\item $\min\{ \diam(\gamma'[f(x),f(u)]), \diam(\gamma'[f(y),f(u)])\}\leq \lambda_1\delta_{G'}(f(u));$
\item $\diam (\gamma')\leq \lambda_2 |f(x)-f(y)|,$\end{enumerate}
where $\gamma'=f(\gamma)$ and $\lambda_i$ depends only on $c,H$ and $b$ for $i=1,2$.
\end{thm}

We remark that with the extra local compactness assumption, it is not hard to see that $G$ is a proper geodesic space with respect to the quasihyperbolic metric, see \cite[Proposition $2.8$]{BHK}, because a complete locally compact length space is proper and geodesic. Hence the above result holds also for quasihyperbolic geodesic of $G$ in this situation.

%

The rest of this paper is organized as follows.  In Section \ref{sec-2}, we recall some definitions and preliminary results, particularly, some basic properties of short arcs. In Section \ref{sec-4}, Theorem \ref{thm1.1} is proved based on the properties of short arcs.
Section \ref{sec-5} is devoted to the proof of Theorem \ref{thm-3}.

\section{Preliminaries}\label{sec-2}

In this section, we give the necessary definitions and auxiliary results, which will be used in the proofs of our main results.

Throughout this paper,
balls and spheres in metric spaces $X$ are written as
$$\mathbb{B}(a,r)=\{x\in X:\,|x-a|<r\},\;\;\mathbb{S}(a,r)=\{x\in X:\,|x-a|=r\}$$
and $$
\mathbb{\overline{B}}(a,r)=\mathbb{B}(a,r)\cup \mathbb{S}(a,r)= \{x\in X:\,|x-a|\leq r\}.$$

For convenience, given
domains $G \subset X,$   $G' \subset Y$, a map $f:G \to G'$
and points $x$, $y$,
$z$, $\ldots$ in  $G$, we always  denote by $x'$, $y'$, $z'$, $\ldots$
the images in $G'$ of $x$, $y$, $z$, $\ldots$ under $f$,
respectively. Also, we assume that $\gamma$
denotes an arc in $G$ and $\gamma'$
 the image in $G'$ of $\gamma$
under $f$.

\subsection{Quasihyperbolic metric, solid arcs and  short arcs}
In this subsection, we start with the definition of quasihyperbolic metric. If $X$ is a connected metric space and $G\varsubsetneq X$ is a non-empty open set,
then it follows from \cite[Rem. 2.2]{HL} that the boundary of $G$ satisfies $\partial G\not=\emptyset$.
Suppose $\gamma\subset G$ denotes a rectifiable arc or a path, its {\it quasihyperbolic length} is the number:

$$\ell_{k_G}(\gamma)=\int_{\gamma}\frac{|dz|}{\delta_{G}(z)},
$$ where $\delta_G(z)$ denotes the distance from $z$ to $\partial G$.


For each pair of points $x$, $y$ in $G$, the {\it quasihyperbolic distance}
$k_G(x,y)$ between $x$ and $y$ is defined in the following way:
$$k_G(x,y)=\inf\ell_{k_G}(\gamma),
$$
where the infimum is taken over all rectifiable arcs $\gamma$
joining $x$ to $y$ in $G$.

 If $\gamma$ is a rectifiable curve in $G$ connecting $x$ and $y$, then (see, e.g., the proof of Theorem $2.7$ in \cite{HL} or \cite{Vai6-0})
\beq\label{base-eq-1}\ell_{k_G}(\gamma)\geq \log\Big(1+\frac{\ell(\gamma)}
{\min\{\delta_{G}(x), \delta_{G}(y)\}}\Big)
\eeq
and thus,
\beq\label{base-eq-2}k_G(x,y)\geq \log\Big(1+\frac{|x-y|}
{\min\{\delta_{G}(x), \delta_{G}(y)\}}\Big).
\eeq

Gehring and Palka \cite{GP} introduced the quasihyperbolic metric of
a domain in $\IR^n$. For the basic properties of this metric we refer to \cite{Geo}. Recall that a curve $\gamma$ from $x$ to
$y$ is a {\it quasihyperbolic geodesic} if
$\ell_{k_G}(\gamma)=k_G(x,y)$. Each subcurve of a quasihyperbolic
geodesic is obviously a quasihyperbolic geodesic. It is known that a
quasihyperbolic geodesic between any two points in a Banach space $X$ exists if the
dimension of $X$ is finite, see \cite[Lem. 1]{Geo}. This is not
true in arbitrary metric spaces (cf. \cite[Ex. 2.9]{Vai6-0}).

Let us recall a result which is useful for our later discussions.

\begin{Lem}\label{ll-11}(\cite[Lem. 2.4]{HWZ}) Let $X$ be a $c$-quasiconvex metric space and let $G\subsetneq X$ be a domain. Suppose that $x$, $y\in G$ and either $|x-y|\leq \frac{1}{3c}\delta_G(x)$ or $k_G(x,y)\leq 1$. Then
\be\label{vvm-2} \frac{1}{2}\frac{|x-y|}{\delta_G(x)}< k_G(x,y) < 3c\frac{|x-y|}{\delta_G(x)}.\ee
\end{Lem}

\noindent Here, we say that $X$ is {\it c-quasiconvex} $(c\geq 1)$ if each pair of points $x$, $y\in X$ can be joined by an arc $\gamma$ in $X$ with length ${\ell}(\gamma)\leq c|x-y|$.

\bdefe \label{def1.4}
 Suppose $\gamma$ is an arc in a domain $G\varsubsetneq X$ and $X$ is a rectifiably connected metric space. The arc may be closed, open or half open. Let $\overline{x}=(x_0,$ $\ldots,$ $x_n)$,
$n\geq 1$, be a finite sequence of successive points of $\gamma$.
For $h\geq 0$, we say that $\overline{x}$ is {\it $h$-coarse} if
$k_G(x_{j-1}, x_j)\geq h$ for all $1\leq j\leq n$. Let $\Phi_k(\gamma,h)$
denote the family of all $h$-coarse sequences of $\gamma$. Set

$$s_k(\overline{x})=\sum^{n}_{j=1}k_G(x_{j-1}, x_j)$$ and
$$\ell_{k_G}(\gamma, h)=\sup \{s_k(\overline{x}): \overline{x}\in \Phi_k(\gamma,h)\}$$
with the agreement that $\ell_{k_G}(\gamma, h)=0$ if
$\Phi_k(\gamma,h)=\emptyset$. Then the number $\ell_{k_G}(\gamma, h)$ is the
{\it $h$-coarse quasihyperbolic length} of $\gamma$.  \edefe

If $X$ is $c$-quasiconvex, then $\ell_{k_G}(\gamma, 0)=\ell_{k_G}(\gamma)$ (see, e.g., \cite[Prop. A.7, Rem. A.13]{BHK} and \cite[Lem. 2.5]{HWZ} ).

\bdefe \label{def1.5} Let $G$ be a proper domain in a rectifiably connected metric space $X$. An arc $\gamma\subset G$
is {\it $(\nu, h)$-solid} with $\nu\geq 1$ and $h\geq 0$ if
$$\ell_{k_G}(\gamma[x,y], h)\leq \nu\;k_G(x,y)$$ for all $x$, $y\in \gamma$.

An arc $\gamma\subset G$ with endpoints $x$ and $y$ is said to be $\varepsilon$-short ($\varepsilon\geq 0$) if $$\ell_{k_G}(\gamma)\leq k_G(x,y)+\varepsilon.$$

Obviously, by the definition of $k_G$, we know that for every $\varepsilon>0$, there exists an arc $\gamma \subset G$ such that $\gamma$ is $\varepsilon$-short, and it is easy to see that  every subarc of an $\varepsilon$-short arc is also $\varepsilon$-short.\edefe

\br
For any pair of points $x$ and $y$ in a proper domain $G$ of Banach space $E$, if the dimension of $E$ is finite, then there exists a quasihyperbolic geodesic in $G$ connecting $x$ and $y$ (see \cite[Lem. 1]{Geo}).  But if the dimension of $E$ is infinite, this property is no longer valid (see, e.g., \cite[Ex. 2.9]{Vai6-0}). In order to overcome this shortcoming in Banach spaces, V\"ais\"al\"a proved the existence of neargeodesics or quasigeodesics (see \cite{Vai6}), and  every quasihyperbolic geodesic is a quasigeodesic. See also \cite{RT}. In metric spaces, we do not know if this existence property is true or not. However, this existence property plays a very important role in the related discussions.
In order to overcome this disadvantage, in this paper, we will exploit the substitution of ``quasigeodesics" replaced by ``short arcs". The class of short arcs was introduced when V\"ais\"al\"a studied properties of Gromov hyperbolic spaces \cite{Vai9} (see also \cite{BH,Herron}), and we see that the existence of such class of arcs is obvious in metric spaces.
\er


By a slight modification of the method used in the proof of \cite[Lem. 6.21]{Vai6}, we get the following  result.

\begin{lem}\label{ll-14} Suppose that $X$ is a $c$-quasiconvex metric space and that $G\varsubsetneq X$ is a domain, and that $\gamma$ is a  $(\nu,h)$-solid arc in $G$ with endpoints $x$, $y$ such that $\min\{\delta_G(x),\delta_G(y)\}=r\geq 3c|x-y|$. Then there is a constant $\mu_1=\mu_1(c,\nu)$ such that $$\diam(\gamma)\leq \max\{\mu_1|x-y|, 2r(e^h-1)\},$$ where ``$\diam$" means ``{\rm diameter}".
\end{lem}

\bpf Without loss of generality, we assume that $\delta_G(y)\geq \delta_G(x)=r$. Denoting $t=|x-y|$ and applying Lemma \Ref{ll-11}, we get $$k_G(x,y)\leq 3ct/r.$$
Let $u\in \gamma$. To prove this lemma, it suffices to show that there exists a constant $\mu_1=\mu_1(c,\nu)$ such that
\be\label{neq-eq-1}|u-x|\leq \max\big\{\frac{\mu_1}{2}|x-y|,r(e^h-1)\big\}.\ee

To this end, we consider two cases. The first case is:  $k_G(u,x)\leq h$. Under this assumption, we see from \eqref{base-eq-2} that
\be\label{dw-2}|u-x|\leq (e^{k_G(u,x)}-1)\delta_G(x)\leq r(e^h-1).\ee

For the remaining case: $k_G(u,x)> h$, we choose a sequence of successive points of $\gamma$: $x=x_0$, $\ldots$, $x_{n}=u$ such that $$ k_G(x_{j-1},x_j)=h \;\;\; {\rm for}\;\;\; j\in\{1,\ldots, n-1\}$$ and $$0< k_G(x_{n-1},x_n)\leq h.$$ Then $n\geq 2$ and
$$(n-1)h\leq \sum_{i=1}^{n-1}k_G(x_{j-1},x_j)\leq \ell_{k_G}(\gamma,h)\leq \nu k_G(x,y)\leq 3c\nu t/r,$$
which shows that $$k_G(x,u)\leq \sum_{i=1}^{n}k_G(x_{j-1},x_j)\leq nh\leq 6c\nu t/r.$$
Let $s=\frac{t}{r}$. Then $s\leq \frac{1}{3c}$ and
$$\frac{|u-x|}{t}\leq \frac{e^{6c\nu s}-1}{s}.$$
Obviously,
the function $g(s)=\frac{1}{s}\big(e^{6c\nu s}-1\big)$ is increasing in
$(0,\infty]$ and $\lim_{s\to 0}{\frac{e^{6c\nu s}-1}{s}}=6c\nu $. Letting $$\mu_1=6c(e^{2\nu}-1)$$ gives
\be\label{dw-1} |u-x|\leq \frac{1}{2}\mu_1t.\ee
It follows from \eqref{dw-2} and \eqref{dw-1} that \eqref{neq-eq-1} holds, and hence the proof of the lemma is complete.
\epf

\begin{lem} \label{mon-4}  Suppose that $X$ is a $c$-quasiconvex metric space and $G\subsetneq X$ is a domain. Suppose, further, that for $x$,
$y\in G$, \begin{enumerate}
\item
$\gamma$ is an $\varepsilon$-short arc in $G$ connecting $x$ and $y$ with $0<\varepsilon\leq \frac{1}{2}k_{G}(x,y)$, and
\item
$|x-y|\leq \frac{1}{3c} \min\{\delta_{G}(x), \delta_{G}(y)\}$. \end{enumerate}
Then $$\ell(\gamma)\leq \frac{9}{2}ce^{\frac{3}{2}}|x-y|.$$\end{lem}
\bpf Without loss of generality, we assume that $\min\{\delta_{G}(x), \delta_{G}(y)\}=\delta_{G}(x)$. It follows from \eqref{base-eq-1} and Lemma \Ref{ll-11} that
$$\log\left(1+\frac{\ell(\gamma)}{\delta_{G}(x)}\right)\leq \ell_{k_{G}}(\gamma)\leq k_{G}(x,y)+\varepsilon \leq  \frac{3}{2}k_{G}(x,y)\leq \frac{9c}{2}\frac{|x-y|}{\delta_{G}(x)}\leq \frac{3}{2}.
$$ Hence,
$$\frac{\ell(\gamma)}{\delta_{G}(x)}\leq e^{\frac{3}{2}}-1.$$
Let $f(t)=t-e^{\frac{3}{2}}\log(1+t)$. Then $f(t)$ is decreasing for $t\in [-1,e^{\frac{3}{2}}-1]$.  In particular, we have $f(t)\leq f(0)=0$  which leads to
$$\frac{\ell(\gamma)}{\delta_{G}(x)}\leq e^{\frac{3}{2}}\log\left(1+\frac{\ell(\gamma)}{\delta_{G}(x)}\right)\leq \frac{9}{2}ce^{\frac{3}{2}}\frac{|x-y|}{\delta_{G}(x)}.$$ Therefore, $$\ell(\gamma)\leq \frac{9}{2}ce^{\frac{3}{2}}|x-y|,$$ as required.\epf

\subsection{Properties of uniform domains}

Let us recall the following useful property of uniform domains.

\begin{Lem}\label{BHK-lem}$($\cite[Lem. 3.12]{BHK}$)$ Suppose $G\subsetneq X$ is a $b$-uniform domain in a rectifiably connected metric space $X$. Then for any $x, y\in G$, we have $$k_G(x,y)\leq 4b^2\log\left(1+\frac{|x-y|}{\min\{\delta_G(x),\delta_G(y)\}}\right).$$
\end{Lem}

We note that Gehring and Osgood \cite{Geo}  characterized uniform domains in terms
of an upper bound for the quasihyperbolic metric in the case of domains in $ {\mathbb R}^n \,$ as follows: a
domain $G$ is {\em uniform} if and only if there exists a constant $C\ge
1$ such that
$$
k_G(x,y)\le C \log\left(1+\frac{|x-y|}{\min\{\delta_G(x),\delta_G(y)\}}\right)
$$
for all $x,y\in G$. As a matter of fact, the above inequality
appeared in \cite{Geo} in a form with an additive constant on the
right hand side: it was shown by Vuorinen \cite[2.50]{Vu2} that the
additive constant can be chosen to be $0$.

The following are the analogues of Lemmas $6.10$ and $6.11$ in \cite{Vai6} in the setting of metric spaces. The proofs are similar.

\begin{Lem}\label{ll-12} Suppose that $G\subsetneq X$ is a $b$-uniform domain in a rectifiably connected metric space $X$, and that $\gamma$ is an arc in $\{x\in G: \delta_G(x)\leq r\}$. If $\gamma$ is $(\nu,h)$-solid, then  $$\diam(\gamma)\leq M_1 r,$$ where $M_1=M_1(b, \nu,h)$.
\end{Lem}


\begin{Lem}\label{ll-13} For all $b\geq 1$, $\nu\geq 1$ and $h\geq 0$, there are constants $0<q_0=q_0(b,\nu,h)<1$ and $M_2=M_2(b,\nu,h)\geq 1$
with the following property:
Suppose that $G$ is a $b$-uniform domain and $\gamma$ is a $(\nu,h)$-solid arc starting at $x_0\in G$. If $\gamma$ contains a point $u$ with $\delta_G(u)\leq q_0\delta_G(x_0)$, then $$\diam(\gamma_{u})\leq M_2\delta_G(u),$$ where $\gamma_{u}=\gamma\setminus \gamma[x_0,u)$.
\end{Lem}


Now, we are ready to prove an analogue of Lemma \ref{ll-14} for uniform domains.

\begin{lem}\label{ll-15} Suppose that $X$ is a $c$-quasiconvex metric space and that $G\varsubsetneq X$ is a $b$-uniform domain, and that $\gamma$ is a $(\nu,h)$-solid arc in $G$ with endpoints $x$, $y$. Let $\delta_G(x_0)=\max_{p\in \gamma}\delta_G(p)$. Then there exist constants
$\mu_2=\mu_2( b, \nu, h)\geq 1$ and $\mu_3=\mu_3( b, c, \nu, h)\geq 1$ such that
\begin{enumerate} \item\label{ma-0-1}
$\diam(\gamma[x,u])\leq \mu_2 \delta_G(u)$ for $u\in \gamma[x,x_0],$ and $\diam(\gamma[y,v])\leq \mu_2 \delta_G(v)$ for $v\in \gamma[y, x_0]$;
\item\label{ma-0-2} $\diam(\gamma)\leq \max\big\{\mu_3 |x-y|, 2(e^h-1)\min\{\delta_G(x),\delta_G(y)\}\big\}.$
\end{enumerate}
 \end{lem}

\bpf We first prove \eqref{ma-0-1}. Obviously, it suffices to prove the first inequality in \eqref{ma-0-1} because the proof for the second one is similar. Let $$\mu_2=\max\Big\{\frac{M_1}{q_0}, M_2\Big\},$$ where $M_1=M_1(b,\nu,h)$ is the constant from Lemma
\Ref{ll-12}, $q_0=q_0(b,\nu,h)$ and $M_2=M_2(b,\nu,h)$ are the constants from Lemma \Ref{ll-13}.

For $u\in \gamma[x,x_0]$, we divide the proof into two cases. If $\delta_G(u)\leq q_0\delta_G(x_0)$, then Lemma \Ref{ll-13} leads to
\be\label{ma-1} \diam(\gamma[x,u])\leq M_2\delta_G(u).\ee If  $\delta_G(u)> q_0\delta_G(x_0)$, then applying Lemma \Ref{ll-12} with the substitution
$r$ replaced by $\delta_G(x_0)$ and $\gamma$ replaced by $\gamma[x,u]$, we easily get
\be\label{ma-2}\diam(\gamma[x,u])\leq M_1\delta_G(x_0)< \frac{M_1}{q_{0}}\delta_G(u).\ee
It follows from \eqref{ma-1} and \eqref{ma-2} that the first assertion in \eqref{ma-0-1} holds, and thus the proof of \eqref{ma-0-1} is complete.

To prove \eqref{ma-0-2}, without loss of generality, we assume that $$\min\{\delta_G(x),\; \delta_G(y)\}= \delta_G(x)\;\;\; {\rm and}\;\;\; \diam(\gamma)> |x-y|.$$
Let
$$\mu_3=\frac{3}{4}\big[1+2(1+6c)(e^{h+4b^2\nu\log(1+4\mu_2)}-1)\big].$$

If $\delta_G(x)\geq 3c|x-y|$, then $(2)$ follows from Lemma \ref{ll-14} since the constant $\mu_1$ in Lemma \ref{ll-14} satisfies $\mu_1< \mu_3$.
Hence, in the following, we assume that \beq\label{add-eq-1}\delta_G(x)<3c|x-y|.\eeq Let $x_1\in \gamma$ (resp. $y_1\in \gamma$) be the first point in $\gamma$ from $x$ to $y$ (resp. from $y$ to $x$) such that (see Figure $1$)
\beq\label{eq-new-1}\diam(\gamma[x,x_1])=\frac{1}{2}|x-y|\;\,(\mbox{resp.}\; \diam(\gamma[y,y_1])=\frac{1}{2}|x-y|).\eeq

%


\noindent Then we have $$\diam(\gamma[y,x_1])\geq |y-x_1|\geq |y-x|-|x-x_1|\geq\frac{1}{2}|y-x|=\diam(\gamma[x,x_1]),$$ and similarly, we get $$\diam(\gamma[x,y_1])>\diam(\gamma[y,y_1]).$$
Thus, it follows from  \eqref{ma-0-1} that $$\frac{1}{2}|x-y|=\diam(\gamma[x,x_1])=\diam(\gamma[y,y_1])\leq  \mu_2
\min\{\delta_G(x_1), \delta_G(y_1)\}.$$ Also,
$$|x_1-y_1|\leq |x_1-x|+|x-y|+|y-y_1|\leq 2|x-y|.$$ Then Lemma \Ref{BHK-lem} implies
$$k_G(x_1,y_1)\leq 4b^2\log\left(1+\frac{|x_1-y_1|}{\min\{\delta_G(x_1),\delta_G(y_1)\}}\right)\leq 4b^2\log(1+4\mu_2).$$
Since $\gamma$ is a $(\nu,h)$-solid arc, for any $u_1,$ $u_2\in \gamma[x_1,y_1]$, we have
\begin{eqnarray*}k_G(u_1,u_2)&\leq&   \max\{h,\; \ell_{k_G}(\gamma[x_1,y_1],h)\}\leq h+\nu k_G(x_1,y_1)\\&\leq& h+4b^2\nu\log(1+4\mu_2),\end{eqnarray*}
and so, for all $z\in \gamma[x_1,y_1]$, we get from \eqref{base-eq-2}, \eqref{add-eq-1} and \eqref{eq-new-1} that
\begin{eqnarray}\label{eq-new-2}|z-x_1|&\leq& (e^{k_G(z,x_1)}-1)\delta_G(x_1)\\\nonumber&\leq&   (e^{h+4b^2\nu\log(1+4\mu_2)}-1)(\delta_G(x)+|x-x_1|)\\\nonumber&\leq& \frac{1}{2}(1+6c)(e^{h+4b^2\nu\log(1+4\mu_2)}-1)|x-y|.\end{eqnarray}
Let $w_1,$ $w_2\in \gamma$ be points such that \be \label{sat-1}|w_1-w_2|\geq \frac{2}{3}\diam(\gamma).\ee
Then we get
\bcl\label{sat-2}
$|w_1-w_2|\leq \frac{2}{3}\mu_3|x-y|$.
\ecl
Since \eqref{eq-new-1} guarantees that neither $\gamma[x,x_1]$ nor $\gamma[y, y_1]$ contains the set $\{w_1,w_2\}$, we see that,
to prove this claim, according to the positions of $w_1$ and $w_2$ in $\gamma$, we need to consider the following four possibilities.\begin{enumerate}
\item
 $w_1\in \gamma[x,x_1]$ and $w_2\in \gamma[y,y_1]$.
 Obviously, by \eqref{eq-new-1}, we have  $$|w_1-w_2|\leq |w_1-x|+|x-y|+|y-w_2|\leq 2|x-y|.$$

\item
 $w_1\in \gamma[x,x_1]$ and $w_2\in \gamma[x_1,y_1]$. Then \eqref{eq-new-1} and \eqref{eq-new-2} show that $$|w_1-w_2|\leq |w_1-x_1|+|x_1-w_2|\leq \frac{1}{2}\big[1+(1+6c)(e^{h+4b^2\nu\log(1+4\mu_2)}-1)\big]|x-y|.$$

\item
$w_1,$ $w_2\in \gamma[x_1,y_1]$. Then \eqref{eq-new-2} implies $$|w_1-w_2|\leq |w_1-x_1|+|x_1-w_2|\leq (1+6c)(e^{h+4b^2\nu\log(1+4\mu_2)}-1)|x-y|.$$

\item
$w_1\in \gamma[x_1,y_1]$ and $w_2\in \gamma[y_1,y]$. Again, we infer  from \eqref{eq-new-1} and \eqref{eq-new-2} that $$ |w_1-w_2|\leq |w_1-x_1|+|x_1-y_1|+|y_1-w_2|\leq \frac{1}{2}\big[1+2(1+6c)(e^{h+4b^2\nu\log(1+4\mu_2)}-1)\big]|x-y|.$$
\end{enumerate}
The claim is proved.\medskip

Now, we are ready to finish the proof. It follows from \eqref{sat-1} and Claim \ref{sat-2} that
$$\diam(\gamma)\leq \frac{3}{2}|w_1-w_2|\leq \mu_3|x-y|,$$
which implies that \eqref{ma-0-2} also holds in this case. Hence, the proof of the lemma is complete.
\epf

\subsection{Free quasiconformal mappings and coarsely quasihyperbolic mappings}

The definition of free quasiconformality is as follows.

\bdefe\label{japan-33} Let $G\varsubsetneq X$ and $G'\varsubsetneq Y$ be two domains (open and connected), and let $\varphi:[0,\infty)\to [0,\infty)$ be a homeomorphism with $\varphi(t)\geq t$. We say
that a homeomorphism $f: G\to G'$ is \begin{enumerate}
\item \label{sunday-1}
{\it $\varphi$-semisolid } if $$k_{G'}(f(x),f(y))\leq \varphi(k_G(x,y))$$
for all $x$, $y\in G$;

\item \label{sunday-2} $\varphi$-{\it solid} if both $f$ and ${\it f^{-1}}$ are $\varphi$-semisolid;

\item {\it freely
$\varphi$-quasiconformal} ($\varphi$-FQC in brief) or {\it fully $\varphi$-solid}
if $f$ is
$\varphi$-solid in every subdomain of $G$,\end{enumerate}
where $k_G(x,y)$ denotes the quasihyperbolic distance of $x$ and $y$ in $G$. See Section \ref{sec-2} for the precise definitions of $k_G(x,y)$ and other
notations and concepts in the rest of this section.
\edefe
\bdefe Let $G\varsubsetneq X$ and $G'\varsubsetneq Y$ be two domains. We say
that a homeomorphism $f: G\to G'$ is \begin{enumerate}

\item
 {\it $C$-coarsely $M$-quasihyperbolic}, or briefly
$(M,C)$-CQH, if there are constants $M\geq 1$ and $C\geq 0$ such that for all $x$, $y\in G$,
$$\frac{k_G(x,y)-C}{M}\leq k_{G'}(f(x),f(y))\leq M\;k_G(x,y)+C.$$
%
\item
{\it fully $C$-coarsely $M$-quasihyperbolic}  if there are constants $M\geq 1$ and $C\geq 0$ such that
$f$ is $C$-coarsely $M$-quasihyperbolic in every subdomain of $G$.

\end{enumerate}\edefe

Under coarsely quasihyperbolic mappings, we have the following useful relationship between short arcs and solid arcs.

\begin{lem}\label{ll-001} Suppose that $X$ and $Y$ are rectifiably connected metric spaces, and that $G\varsubsetneq X$ and $G'\varsubsetneq Y$ are domains. If $f:\;G\to G'$ is $(M,C)$-CQH, and $\gamma$ is an $\varepsilon$-short arc in $G$ with $0<\varepsilon\leq 1$, then there are constants $\nu=\nu(C,M)$ and $h=h(C, M)$ such that
 the image $\gamma'$ of $\gamma$ under $f$ is $(\nu,h)$-solid in $G'$.
\end{lem}
\bpf
Let$$h=(2M+1)C+2M\;\; \mbox{and}\;\; \nu=\frac{4(C+1)M(M+1)}{2C+1}.$$
Obviously, we only need to verify that for $x$, $y\in \gamma$,
\be\label{new-eq-3}\ell_{k_{G'}}(\gamma'[x',y'],h)\leq\nu k_{G'}(x',y').\ee We prove this by considering two cases.
The first case is: $k_G(x,y)<2C+1$. Then for $z_1$, $z_2\in\gamma[x, y]$, we have
$$k_{G'}(z'_1,z'_2)\leq Mk_G(z_1,z_2)+C\leq M(k_G(x,y)+\varepsilon)+C<(2M+1)C+2M=h,$$
and so
\be\label{ma-3}\ell_{k_{G'}}(\gamma'[x',y'],h)=0.\ee

Now, we consider the other case: $k_G(x,y)\geq 2C+1$. Then $$k_{G'}(x',y')\geq \frac{1}{M}(k_G(x,y)-C)> \frac{1}{2M}k_G(x,y).$$ With the aid of \cite[Theorems 4.3 and 4.9]{Vai6}, we have
\beq\label{ma-4}
\ell_{k_{G'}}(\gamma'[x',y'],h) &\leq&
\ell_{k_{G'}}(\gamma'[x',y'],(M+1)C) \leq
 (M+1)\ell_{k_G}(\gamma[x,y])\\ \nonumber &\leq&(M+1)(k_G(x,y)+\varepsilon)\\ \nonumber &\leq& \frac{2(C+1)(M+1)}{2C+1}k_{G}(x,y) \\ \nonumber &\leq& \frac{4(C+1)M(M+1)}{2C+1}k_{G'}(x',y').\eeq
It follows from \eqref{ma-3} and \eqref{ma-4} that \eqref{new-eq-3} holds.\epf


The following results will be used in the proof of Theorem \ref{thm1.1}.

\begin{lem}\label{ll-000} Suppose that $X$ and $Y$ are both $c$-quasiconvex and complete metric \;\; spaces, and that $G\varsubsetneq X$ and $G'\varsubsetneq Y$ are domains. If both $f:$ $G\to G'$ and $f^{-1}:$ $G'\to G$ are weakly $H$-quasisymmetric, then
\begin{enumerate}
\item \label{sat-3}
$f$ is $\varphi$-FQC, where $\varphi=\varphi_{c,H}$ which means that the function $\varphi$ depends only on $c$ and $H$;

\item \label{sat-4}
 $f$ is fully $(M,C)$-CQH, where $M=M(c,H)\geq 1$ and $C=C(c,H)\geq 0$ are constants.
 \end{enumerate}
\end{lem}
\bpf By \cite[Thm. 1.6]{HL}, we know that for every subdomain $D\subset G$, both $f:$ $D\to D'$ and $f^{-1}:$ $D'\to D$ are $\varphi$-semisolid with $\varphi=\varphi_{c,H}$, and so, $f$ is $\varphi$-FQC. Hence \eqref{sat-3} holds.  On the other hand,  \cite[Theorem 1]{ HWZ} implies that \eqref{sat-3} and \eqref{sat-4} are equivalent, and thus, \eqref{sat-4} also holds.
\epf

\begin{Lem}\label{lem-ll-0}$($\cite[Lem. 6.5]{Vai8}$)$ Suppose that $X$ is $c$-quasiconvex, and that $f:$ $X\to Y$ is weakly $H$-quasisymmetric. If $x,$ $y,$ $z$ are distinct points in $X$ with $|y-x|\leq t|z-x|$, then $$|y'-x'|\leq \theta(t)|z'-x'|,$$ where the function $\theta(t)=\theta_{c, H}(t)$ is increasing in $t$.
\end{Lem}

\begin{Lem}\label{lem-ll-1}$($\cite[Lem. 5.4]{Vai6-0}$)$ Suppose that $f:$ $X\to Y$ is weakly $H$-quasisymmetric and that $f(X)$ is $c$-quasiconvex. If $x,$ $y,$ $z$ are distinct points in $X$ with $|y-x|= t|z-x|$ and if $0<t\leq 1$, then $$|y'-x'|\leq \theta(t)|z'-x'|,$$ where $\theta:[0,1]\to[0,\infty)$ is an embedding with $\theta(0)=0$ depending only on $H$ and $c$.
\end{Lem}

The following result easily follows from Lemma \Ref{lem-ll-1}.

\begin{lem}\label{lem-ll-2}  Suppose that $f:$ $X\to Y$ is weakly $H$-quasisymmetric and that $f(X)$ is $c$-quasiconvex. Then $f^{-1}:$ $f(X)\to X$ is weakly $H_1$-quasisymmetric with $H_1$ depending only on $H$ and $c$.
\end{lem}
\bpf Let  $x,$ $y,$ $z$ be distinct points in $X$ with $|y-x|= t|z-x|$ and  $0<t\leq 1$. Then by Lemma \Ref{lem-ll-1} there exists some constant $c_1>2$ such that $$|y'-x'|\leq \theta(t)|z'-x'|,$$ where $\theta(\frac{1}{c_1})<1$.

In order to prove that $f^{-1}:$ $f(X)\to X$ is weakly $H$-quasisymmetric, we let $a,b,x$ be distinct points in $X$ with $|a'-x'|\leq |b'-x'|$. Then if $|a-x|\leq c_1|b-x|$, there is nothing to prove. Hence, we assume that $|a-x|> c_1|b-x|$.
 Then we have
$$|a'-x'|\leq |b'-x'|\leq\theta(1/c_1)|a'-x'|<|a'-x'|.$$
This contradiction  completes the proof.
\epf



\section{The proof of Theorem \ref{thm1.1}}\label{sec-4}

In this section, we always assume that $X$ and $Y$ are $c$-quasiconvex and complete metric spaces, and that $G\varsubsetneq X$ and $G'\varsubsetneq Y$ are domains. Furthermore, we suppose that  $f:$ $G\to G'$ is weakly $H$-quasisymmetric, $G'$ is $c_1$-quasiconvex and $D\subset G$ is $b$-uniform.

Under these assumptions, it follows from Lemmas \ref{ll-000} and  \ref{lem-ll-2}   that $f$ is $(M, C)$-$CQH$ with $M=M(c, H)\geq 1$ and $C=C(c, H)\geq 0$.

We are going to show the uniformity of $D'=f(D)$. For this, we let $x'$, $y'\in D'=f(D)\subset G'$, and $\gamma'$ be an $\varepsilon$-short arc in $D'$ joining $x'$ and $y'$ with
$$0<\varepsilon<\min\big\{1,\frac{1}{2}k_{D'}(x',y')\big\}.$$
Then by Lemma \ref{ll-001}, the preimage $\gamma$ of $\gamma'$ is a $(\nu,h)$-solid arc in $D$ with $\nu=\nu(c, H)$ and $h=h(c, H)$. Let $w_0\in\gamma$ be such that (see Figure $2$)
\be\label{wes-1}
\delta_D(w_0)=\max_{p\in \gamma}\delta_D(p).
\ee

%

\noindent Then by Lemma \ref{ll-15}, there is a constant $\mu=\mu(b,\nu,h)$ such that for each $u\in\gamma[x,w_0]$ and for all $z\in\gamma[u,w_0]$, \be\label{new-eq-4}|u-z|\leq \diam (\gamma[u, z])
\leq\mu\delta_D(z),\ee
and for each $v\in\gamma[y,w_0]$ and for all $z\in\gamma[v,w_0]$, \be\nonumber |v-z|\leq \diam (\gamma[v, z])\leq\mu\delta_D(z).\ee

In the following, we show that $\gamma'$ is a double cone arc in $D'$. Precisely, we shall prove
that there exist constants $A\geq 1$ and $B\geq 1$ such that for every $z'\in\gamma'[x',y']$,
\be\label{main-eq-1}\min\{\ell(\gamma'[x',z']),\ell(\gamma'[z',y'])\}\leq A\delta_{D'}(z')\ee
and
\be\label{main-eq-2}\ell(\gamma')\leq B|x'-y'|.\ee

The verification of \eqref{main-eq-1} and \eqref{main-eq-2} is given in the following two subsections.

\subsection{The proof of \eqref{main-eq-1}}

Let
$$A=2e^{8b^2A_1(C+1)M}\;\; {\rm and}\;\; A_1=2e^{M+C}(1+\mu)\theta''\Big(6c\theta'(\mu)e^{4b^2M+C}\Big),$$ where the functions $\theta'=\theta'_{b,H}$ and $\theta''=\theta''_{c_1,H}$ are from Lemma \Ref{lem-ll-0}.  Obviously, we only need to get the following estimate:
for all $z'\in\gamma'[x',w'_0]$ (resp. $z'\in\gamma'[y',w'_0]$),
\be\label{main-eq-3}
\ell(\gamma'[x',z'])\leq A\delta_{D'}(z')\; \; ({\rm resp.}\; \ell(\gamma'[y',z'])\leq A\delta_{D'}(z')).
\ee

It suffices to prove the case $z'\in\gamma'[x',w'_0]$ since the proof of the case $z'\in\gamma'[y',w'_0]$ is similar.
Suppose on the contrary that there exists some point $x'_0\in \gamma'[x',w'_0]$ such that
\be\label{main-eq-4}
\ell(\gamma'[x',x_0'])>A\delta_{D'}(x_0').
\ee Then we choose $x'_1\in\gamma'[x',w'_0]$ to be the first point from $x'$ to $w_0'$ such that (see Figure $3$)
\be\label{sat-5}
\ell(\gamma'[x',x'_1])=A\delta_{D'}(x'_1).
\ee

%


Let $x_2\in  D$ be such that (see Figure $2$)
$$|x_1-x_2|=\frac{1}{2}\delta_{D}(x_1).$$


\noindent Then we have
\bcl\label{eq-lwz4}
$|x_1'-x_2'|< e^{4b^2M+C}\delta_{D'}(x_1').$
\ecl
Obviously,
$$\delta_{D}(x_2)\geq \delta_{D}(x_1)-|x_1-x_2|=|x_1-x_2|,$$
and so, \eqref{base-eq-2} and Lemma \Ref{BHK-lem} imply
\begin{eqnarray*}\log\left(1+\frac{|x_1'-x_2'|}{\delta_{D'}(x_1')}\right)&\leq& k_{D'}(x_1',x_2')\leq Mk_D(x_1,x_2)+C\\&\leq& 4b^2M\log\left(1+\frac{|x_1-x_2|}{\min\{\delta_D(x_1),\delta_D(x_2)\}}\right)+C\\&<& 4b^2M+C,
\end{eqnarray*}
whence
$$
|x_1'-x_2'|< e^{4b^2M+C}\delta_{D'}(x_1'),
$$ which shows that the claim holds.\medskip

Let $x'_3\in\gamma'[x',x'_1]$ be such that
\be\label{sun-2}\ell(\gamma'[x',x'_3])=\frac{1}{2}\ell(\gamma'[x',x'_1]),\ee
and then, we get an estimate on $|x_1'-x_2'|$ in terms of $d_{D'}(x_3')$ as stated in the following claim.
\bcl\label{eq-lwz1}
$|x_1'-x_2'|< 2e^{4b^2M+C}\delta_{D'}(x_3').$
\ecl
It follows from \eqref{sat-5} and \eqref{sun-2} that
$$
\delta_{D'}(x_1')< 2 \delta_{D'}(x_3'),
$$ since the choice of $x_1'$ implies $\ell(\gamma'[x',x'_3])<Ad_{D'}(x_3')$. Hence, Claim \ref{eq-lwz4} leads to
$$
|x_1'-x_2'|< 2e^{4b^2M+C}\delta_{D'}(x_3'),
$$ as required.\medskip

On the basis of  Claim \ref{eq-lwz1}, we have
\bcl\label{sun-1}
$|x'_1-x'_3|\leq 2\theta'(\mu)e^{4b^2M+C}\delta_{D'}(x'_3).$
\ecl

In order to apply Lemma \Ref{lem-ll-0} to prove this claim, we need some preparation. It follows from \eqref{base-eq-1}, \eqref{sat-5} and \eqref{sun-2} that
\begin{eqnarray*} k_{D'}(x'_1,x'_3) &\geq& \ell_{k_{D'}}(\gamma'[x'_1,x'_3])-\varepsilon
\geq \log\Big(1+\frac{\ell(\gamma'[x'_1,x'_3])}{\delta_{D'}(x'_1)}\Big)-1
\\ \nonumber&=&  \log\big(1+\frac{A}{2}\big)-1.
\end{eqnarray*} Hence, by Lemma \Ref{BHK-lem}, we have
\begin{eqnarray*}\log\left(1+\frac{|x_1-x_3|}{\min\{\delta_D(x_1),\delta_D(x_3)\}}\right)&\geq&\frac{1}{4b^2}k_D(x_1,x_3)\geq \frac{1}{4b^2M}(k_{D'}(x'_1,x'_3)-C)\\&\geq&
\frac{1}{4b^2M}\Big(\log\big(1+\frac{A}{2}\big)-1-C\Big)\\ \nonumber
&>&\log(1+A_1),\end{eqnarray*}
and so \be\label{z-006}|x_1-x_3|>A_1\min\{\delta_D(x_1),\delta_D(x_3)\}>\frac{A_1}{1+\mu}\delta_D(x_3),\ee which by \eqref{new-eq-4} implies
$$\delta_D(x_3)\leq \delta_D(x_1)+|x_1-x_3|\leq (1+\mu)\delta_D(x_1).$$
Again, by \eqref{new-eq-4}, we know
$$|x_1-x_3|\leq\mu\delta_D(x_1)= 2\mu|x_1-x_2|.$$

Now, we are ready to apply Lemma \Ref{lem-ll-0} to the points $x_1$, $x_2$ and $x_3$ in $D$. Since $f$ is weakly $H$-quasisymmetric and $D$ is $b$-uniform, by considering the restriction $f|_D$ of $f$ onto $D'$, we know from Lemma \Ref{lem-ll-0} that there is an increasing function $\theta'=\theta'_{b,H}$ such that
$$|x'_1-x'_3|\leq \theta'(2\mu)|x'_1-x'_2|,$$ and thus, Claim
\ref{eq-lwz1} assures that
$$
|x'_1-x'_3|\leq  2\theta'(2\mu)e^{4b^2M+C}\delta_{D'}(x'_3),
$$
which completes the proof of Claim \ref{sun-1}.
\medskip

Let us proceed with the proof. To get a contradiction to the contrary assumption \eqref{main-eq-4}, we choose $x'_4\in  D'$ such that
\be\label{sun-3} |x'_3-x'_4|=\frac{1}{3c}\delta_{D'}(x'_3).\ee
Then Lemma \Ref{ll-11} implies that \begin{eqnarray*}\log\left(1+\frac{|x_3-x_4|}{\delta_{D}(x_3)}\right)&\leq& k_{D}(x_3,x_4)\leq Mk_{D'}(x_3',x_4')+C\\ \nonumber&\leq& 3cM\frac{|x_3'-x_4'|}{\delta_{D'}(x_3')}+C\leq M+C,\end{eqnarray*} which yields that \be\label{sun-3-1}|x_3-x_4|< e^{M+C}\delta_{D}(x_3).\ee
On the other hand, Claim \ref{sun-1} and \eqref{sun-3} imply that  $$|x_1'-x_3'|\leq 2\theta'(2\mu)e^{4b^2M+C}\delta_{D'}(x_3')= 6c\theta'(2\mu)e^{4b^2M+C}|x_3'-x_4'|.$$
Now, we apply Lemma \Ref{lem-ll-0} to the points $x_1'$, $x_3'$ and $x_4'$ in $G'$.
Since by Lemma \ref{lem-ll-2} $f^{-1}:$ $G'\to G$ is weakly $H$-quasisymmetric and $G'$ is $c_1$-quasiconvex,  we know from Lemma \Ref{lem-ll-0} that there is an increasing function $\theta''=\theta''_{c_1,H}$ such that
$$|x_1-x_3|\leq \theta''\Big(6c\theta'(2\mu)e^{4b^2M+C}\Big)|x_3-x_4|,$$
which, together with \eqref{z-006} and \eqref{sun-3-1}, shows that
\begin{eqnarray*}
|x_1-x_3|&\leq&  e^{M+C}\theta''\Big(6c\theta'(2\mu)e^{4b^2M+C}\Big)\delta_D(x_3)\\ \nonumber&\leq& \frac{1+\mu}{A_1}e^{M+C}\theta''\Big(6c\theta'(2\mu)e^{4b^2M+C}\Big)|x_1-x_3|\\ \nonumber&=&
\frac{1}{2}|x_1-x_3|.
\end{eqnarray*}
This obvious contradiction shows that \eqref{main-eq-1} is true.
\qed

\subsection{The proof of \eqref{main-eq-2}}
 Let
$$B=12cA^2e^{6b^2M\mu\big(1+\theta''\big(\frac{1+12cA}{3c}\big)\big)},$$
and suppose on the contrary that \be\label{eq-lwz3}\ell(\gamma')> B|x'-y'|.\ee
Since $\frac{9}{2}ce^{\frac{3}{2}}<B$, we see from Lemma \ref{mon-4} that
\be \label{eq-lwz2}|x'-y'|>\frac{1}{3c} \max\{\delta_{D'}(x'), \delta_{D'}(y')\}.\ee
For convenience, in the following, we assume that $$\max\{\delta_{D'}(x'), \delta_{D'}(y')\}=\delta_{D'}(x').$$

First, we choose some special points from $\gamma'$.
By \eqref{eq-lwz3}, we know that there exist $w'_1$ and $w'_2\in \gamma'$ such that $x'$, $w'_1$, $w'_2$ and $y'$ are successive points in $\gamma'$ and
\be\label{eq-11-1}\ell(\gamma'[x',w'_1])=\ell(\gamma'[w'_2,y'])=6cA|x'-y'|.\ee
Then we have
\bcl\label{zzz-002} $|x'-w'_1|\geq \frac{1}{2}\delta_{D'}(w'_1)$ and $|y'-w'_2|\geq \frac{1}{2}\delta_{D'}(w'_2)$.
\ecl

Obviously, it suffices to show the first inequality in the claim. Suppose
 $$|x'-w'_1|< \frac{1}{2}\delta_{D'}(w'_1).$$ Then \eqref{main-eq-1} and \eqref{eq-lwz2} lead to
$$\delta_{D'}(x')\geq \delta_{D'}(w'_1)-|x'-w'_1|>\frac{1}{2}\delta_{D'}(w'_1)\geq \frac{1}{2A}\ell(\gamma'[x',w'_1])=3c|x'-y'|>\delta_{D'}(x').$$
This obvious contradiction completes the proof of Claim \ref{zzz-002}.\medskip

By using Claim \ref{zzz-002}, we get a lower bound for $|w_1-w_2|$ in terms of $\min\{\delta_D(w_1),\; \delta_D(w_2)\}$, which is as follows.
\bcl\label{eq-l2}
$|w_1-w_2|> \Big(1+\theta''\Big(\frac{1+12cA}{3c}\Big)\Big)\mu\min\{\delta_D(w_1),\; \delta_D(w_2)\}.$
\ecl
Without loss of generality, we assume that $\min\{\delta_D(w_1),\; \delta_D(w_2)\}=\delta_{D}(w_1)$.
Then by \eqref{eq-11-1} and Claim \ref{zzz-002}, we have
\be\label{eq-l10}\delta_{D'}(w'_1)\leq 2|x'-w'_1|\leq 2\ell(\gamma'[x',w'_1])=12cA|x'-y'|.\ee
Since $\gamma'$ is an $\varepsilon$-short arc and $D$ is $b$-uniform, by Lemma \Ref{BHK-lem}, we have
\begin{eqnarray*}
\log\left(1+\frac{|w_1-w_2|}{\delta_D(w_1)}\right) &\geq& \frac{1}{4b^2}k_D(w_1,w_2) \geq \frac{1}{4Mb^2}k_{D'}(w'_1,w'_2)-\frac{C}{4Mb^2}\\&\geq& \frac{1}{4Mb^2}\ell_{k_{D'}}(\gamma'[w'_1,w'_2])-\frac{\varepsilon+C}{4Mb^2}
\\ \nonumber&\geq& \frac{1}{4Mb^2}\log\left(1+\frac{\ell(\gamma'[w'_1,w'_2])}{\delta_{D'}(w'_1)}\right)-\frac{1+C}{4Mb^2}
\\ \nonumber&\geq&  \frac{1}{4Mb^2}\log\Big(1+\frac{B-12cA}{12cA}\Big)-\frac{1+C}{4Mb^2}
\\ \nonumber &=&\lambda,
\end{eqnarray*}
where the last inequality follows from \eqref{eq-l10} and the following inequalities:
$$\ell(\gamma'[w_1',w_2'])=\ell(\gamma')-\ell(\gamma'[x',w_1'])-\ell(\gamma'[y',w_2'])>(B-12cA)|x'-y'|.$$ Hence $$|w_1-w_2|\geq(e^{\lambda}-1)\delta_D(w_1)> \Big(1+\theta''\Big(\frac{1+12cA}{3c}\Big)\Big)\mu\delta_D(w_1),$$ as required. \medskip

Next, we get the following upper bound for $|w_1-w_2|$ in terms of $\min\{\delta_D(w_1),\; \delta_D(w_2)\}$.
\bcl\label{mon-3}
$|w_1-w_2|\leq \theta''\left(\frac{1+12cA}{3c}\right)\mu \min\{\delta_D(w_1),\; \delta_D(w_2)\}.$
\ecl
First, we see that $w_0\in\gamma[w_1,y]$ (see Figure $4$), where $w_0$ is the point in $\gamma$ which satisfies \eqref{wes-1} (see Figure $2$), because otherwise \eqref{new-eq-4} gives that $$|w_1-w_2|\leq\mu\delta_D(w_1),$$ which contradicts Claim \ref{eq-l2}.

%


We are going to apply Lemma \Ref{lem-ll-0} to the points $x'$, $w_1'$ and $w_2'$ in $G'$. We need a relationship between $|w'_1-w'_2|$ and $|x'-w_1'|$.
To this end, it follows from \eqref{eq-11-1} that $$|w'_1-w'_2|\leq |w'_1-x'|+|x'-y'|+|y'-w'_2|\leq (1+12cA)|x'-y'|\leq\frac{1+12cA}{3c}|x'-w_1'|,$$
since we infer from
the choice of $w'_1$, \eqref{main-eq-1} and Claim \ref{zzz-002} that
$$|x'-w_1'|\geq \frac{1}{2}\delta_{D'}(w_1')\geq \frac{1}{2A}\ell(\gamma'[x',w'_1])=3c|x'-y'|.$$
Then by Lemma  \Ref{lem-ll-0}, we know that there is an increasing function $\theta''=\theta''_{c_1,H}$ such that
 $$|w_1-w_2|\leq \theta''\left(\frac{1+12cA}{3c}\right)|x-w_1|,$$
and thus, \eqref{new-eq-4} leads to
$$|w_1-w_2|\leq  \theta''\left(\frac{1+12cA}{3c}\right)\mu \min\{\delta_D(w_1),\; \delta_D(w_2)\},$$
which shows that Claim \ref{mon-3} holds.\medskip

We observe that Claim \ref{eq-l2}  contradicts Claim \ref{mon-3}, which completes the proof of \eqref{main-eq-2}.

Inequalities \eqref{main-eq-1} and \eqref{main-eq-2}, together with the arbitrariness of the choice of $x'$ and $y'$ in $D'$, show that $D'$ is $B$-uniform, which implies that Theorem \ref{thm1.1} holds.\qed

\bigskip

\section{The proof  of Theorem  \ref{thm-3}}\label{sec-5}

First, one easily sees from Lemma \ref{lem-ll-2} that $f^{-1}$ is also weakly $H_1$-quasisymmetric for some constant $H_1\geq 1$ because every $b$-uniform domain is clearly $b$-quasiconvex. Thus by Lemma \ref{ll-000}, we get that $f$ is $(M,C)$-CQH for some $M,C\geq 1$. Moreover, using Lemma \ref{ll-001}, we obtain that the image curve $\gamma'$ is $(\nu,h)$-solid for some $\nu,h\geq 1$. Hence the first assertion follows from Lemma \ref{ll-15}.

So, it remains to show the second assertion. We note that again by Lemma \ref{ll-15} we have
\beq\label{eq-li-new1}\diam (\gamma')\leq \mu_3 \max\{ |x'-y'|,  \min\{\delta_{G'}(x'),  \delta_{G'}(y')\} \}.\eeq
 Without loss of generality, we may assume $\delta_{G'}(x')\leq \delta_{G'}(y')$. Since $f^{-1}$ is also weakly $H_1$-quasisymmetric, by \cite[(3.10)]{HL}, one immediately sees that there is an increasing continuous function $\theta:(0,\frac{1}{54c})\to \mathbb{R}_+$ with $\theta$ depending only on $c$ and $H_1$ and with $\theta(0)=0$ such that
\beq\label{eq-li-new2}\frac{|x-y|}{\delta_G(x)}\leq \theta\Big(\frac{|x'-y'|}{\delta_{G'}(x')}\Big),\eeq
whenever $x',y'\in G'$ with $|x'-y'|\leq \frac{\delta_{G'}(x')}{54c}$. Thus, there is a constant $\lambda\in (0,\frac{1}{54c})$ such that $\theta(\lambda)<\frac{1}{50c^2}$.

Thus we divide the proof of the second assertion into two cases.

 \bca\label{ca1}$|x'-y'|\geq \lambda \delta_{G'}(x')$.\eca

 Then, obviously, we get from \eqref{eq-li-new1} that
$$\diam (\gamma')\leq \frac{\mu_3}{\lambda}  |x'-y'|,$$
as desired.

    \bca\label{ca2}  $|x'-y'|< \lambda \delta_{G'}(x')$.\eca

     Then \eqref{eq-li-new2} and the choice of $\lambda$ imply that $$|x-y|\leq \theta(\lambda)\delta_G(x)<\frac{\delta_G(x)}{50c^2}.$$ We claim that
\beq\label{eq-li-new3}\gamma\subset \mathbb{B}(x,8c|x-y|)\subset \mathbb{B}(x,\frac{\delta_G(x)}{6c}).\eeq
Otherwise, there is some point $w\in\gamma$ such that $|w-x|=8c|x-y|\leq \frac{\delta_G(x)}{6c}$ and $\gamma[x,w]\subset  \mathbb{\overline{B}}(x,8c|x-y|)$. Moreover, since $\gamma$ is a $\varepsilon$-short arc, by Lemma \Ref{ll-11} we have
$$\ell_{k_G}(\gamma)\leq k_G(x,y)+\varepsilon\leq \frac{7}{6}k_G(x,y)\leq \frac{7c|x-y|}{2\delta_G(x)},$$
and
$$\ell_{k_G}(\gamma)\geq k_G(x,w)\geq \frac{1}{2}\frac{|x-w|}{\delta_G(x)}\geq 4c \frac{|x-y|}{\delta_G(x)},$$
which is an obvious contradiction. Hence the required relation \eqref{eq-li-new3} follows.

Moreover, by \eqref{eq-li-new3}  we have for all $z\in\gamma$, $$|x-z|\leq 8c|x-y|.$$ Since $X$ is $c$-quasiconvex, there is a curve $\beta$ joining $x$ and $z$ with $\ell(\beta)\leq c|x-z|$. Indeed, we know that $\beta\subset G$, because
$$\ell(\beta)\leq c|x-z|\leq 8c^2|x-y|\leq \frac{\delta_G(x)}{6}.$$
Then we shall define inductively the successive points $x=x_0,...,x_n=z$ such that each $x_i$ denotes the last point of $\beta$ in $\overline{\mathbb{B}}(x_{i-1},|x-y|)$, for $i\geq 1$. Obviously, $n\geq 2$, and
$$|x_{i-1}-x_{i}|=|x-y|\;\mbox{ for}\; 1\leq i\leq n-1,\;\;   \;\; |x_{n-1}-x_{n}|\leq |x-y|.$$
Next, we are going  to obtain an upper bound for $n$.
 Since for $1\leq i \leq n-1$, $$\ell(\beta[x_{i-1},x_{i}])\geq |x_{i-1}-x_{i}|=|x-y|,$$ we have
$$(n-1)|x-y|\leq \ell(\beta)\leq c|x-z|\leq 8c^2|x-y|,$$
which implies that
$$n\leq 8c^2+1.$$

Now, we are going to complete the proof in this case. Since for $1\leq i\leq n-1$, $|x_{i+1}-x_{i}|\leq |x_{i-1}-x_{i}|$ and $f$ is weakly $H$-quasisymmetric, we have
$$|x'_1-x'_0|\leq H|x'-y'|,$$
$$|x'_2-x'_1|\leq H|x'_1-x'_0|\leq  H^2|x'-y'|,$$
$$...$$
$$|x'_n-x'_{n-1}|\leq H^n|x'-y'|.$$
Thus,
\begin{eqnarray*}
|x'-z'|&=&|x'_0-x'_n| \leq \sum_{i=1}^n |x'_i-x'_{i-1}|
\\ \nonumber &\leq& (H+H^2+\cdot\cdot\cdot+H^n)|x'-y'|
\\ \nonumber &\leq&  n H^n|x'-y'|
\\ \nonumber &\leq&  (8c^2+1)H^{8c^2+1}|x'-y'|.
\end{eqnarray*}
Therefore, we obtain
$$\diam(\gamma')\leq 2(8c^2+1)H^{8c^2+1}|x'-y'|$$ as desired.

Let $\lambda_2=\max\{\frac{\mu_3}{\lambda}, 2(8c^2+1)H^{8c^2+1}\}$.
Then the proof of this theorem is complete.\\

{\bf Acknowledgement.} The authors would like to thank Professor Xiantao Wang and Professor Manzi Huang for several comments on this manuscripts and thank the referee who has made valuable comments on this manuscripts.

\end{document}